\documentclass[12pt,reqno]{amsart}
\usepackage{amssymb,delarray}
\usepackage{amsfonts,amscd}
\usepackage{epsfig}
\usepackage[all]{xy}

\def\leq{\leqslant}
\def\geq{\geqslant}

\newtheorem{thm}{Theorem}
{Lemma}
\newtheorem{prop}
{Proposition}
\newtheorem{claim}
{Claim}

\newtheorem{cor}
{Corollary}
{Remark}
{Question}
\newtheorem{ex-thm}{Theorem-Example}

{\catcode`\@=11
\gdef\n@te#1#2{\leavevmode\vadjust{%
 {\setbox\z@\hbox to\z@{\strut#1}%
  \setbox\z@\hbox{\raise\dp\strutbox\box\z@}\ht\z@=\z@\dp\z@=\z@%
  #2\box\z@}}}
\gdef\leftnote#1{\n@te{\hss#1\quad}{}}
\gdef\rightnote#1{\n@te{\quad\kern-\leftskip#1\hss}{\moveright\hsize}}
\gdef\?{\FN@\qumark}
\gdef\qumark{\ifx\next"\DN@"##1"{\leftnote{\rm##1}}\else
 \DN@{\leftnote{\rm??}}\fi{\rm??}\next@}}

\begin{document}
\baselineskip=13.7pt plus 2pt 

\title[On the  almost generic covers of the projective plane]{On the almost generic covers of the projective plane}

\author[Vik.S. Kulikov]{Vik.S. Kulikov}

\address{Steklov Mathematical Institute of Russian Academy of Sciences, Moscow, Russia}
 \email{kulikov@mi.ras.ru}

\dedicatory{} \subjclass{}


\maketitle

\def\st{{\sf st}}

\quad \qquad \qquad
{\em Dedicated to G.-M. Greuel on the occasion of his 75th birthday}

\begin{abstract}
A finite morphism $f:X\to \mathbb P^2$ of a a smooth irreducible projective surface $X$ is called an almost generic cover if for each point $p\in \mathbb P^2$ the fibre $f^{-1}(p)$ is supported at least on $\deg f-2$ distinct points
and $f$ is ramified with multiplicity two at a generic point of its ramification locus $R$. In the article, the singular points of the branch curve $B\subset\mathbb P^2$ of an almost generic cover are investigated and main invariants of the covering surface $X$ are calculated in terms of invariants of the curve $B$.
\end{abstract} \vspace{0.5cm}

\def\st{{\sf st}}


\setcounter{section}{-1}

\section{Introduction}
Let $X$ be a smooth irreducible projective surface. A finite morphism $f:X\to \mathbb P^2$, branched along a curve $B\subset \mathbb P^2$,  is called a {\it generic cover of the projective plane} if it has the following properties:
\begin{itemize}
\item[($G_1$)] for each point $p\in \mathbb P^2$
the fibre $f^{-1}(p)$ is supported on at least $\deg f-2$ distinct points,
\item[($G_2$)] $f$ is ramified with multiplicity $2$ at a generic point of its ramification locus $R$,
\item[($G_3$)] the singular points of $B$ are only the ordinary nodes and ordinary cusps.
\end{itemize}

If $X$ is imbedded in some projective space $\mathbb P^n$ then it is well known (see, for example, \cite{C-F}) that the restriction $f:X\to \mathbb P^2$ to $X$ of a linear projection $\text{pr}: \mathbb P^n\to\mathbb P^2$ generic with respect to the imbedding of $X$, is a generic cover of $\mathbb P^2$. Properties of generic covers of $\mathbb P^2$
were investigated  in \cite{K1} -- \cite {K3} and \cite{N}. In this article we generalize the notion of generic covers of
$\mathbb P^2$ as follows. We say that a finite morphism $f:X\to \mathbb P^2$
is an {\it almost generic cover} of the projective plane if it satisfies properties $(G_1)$ and $(G_2)$ of generic covers.

 A dominant morphism $f:X\to\mathbb P^2$ defines a homomorphism $f_*:\pi_1(\mathbb P^2\setminus B,p)\to \mathbb S_{\deg f}$ (called the {\it monodromy} of $f$) whose image $G_f:=f_*(\pi_1(\mathbb P^2\setminus B,p))$  is the {\it monodromy group} of  $f$ and it is a subgroup of the symmetric group $\mathbb S_{\deg f}$ acting on the fibre $f^{-1}(p)=\{ q_1,\dots,q_{\deg f}\}$.

Let $o$ be a point of a curve $B\subset \mathbb P^2$. It is well known that the group $\pi_1^{loc}(B,o):=\pi_1(V\setminus B)$ does not depend on $V$, where $V\subset \mathbb P^2$ is a sufficiently small complex analytic neighbourhood
biholomorphic to a ball of small radius centered at $o$. The image $G_{f,o}:=\text{im}\, f_*\circ i_*$ is called the {\it local monodromy group} of $f$ at the point $o$, where $i_*:\pi_1^{loc}(B,o)=\pi_1(V\setminus B,q)\to \pi_1(\mathbb P^2\setminus B,p)$ is a homomorphism defined (uniquely up to conjugation) by the imbedding $V\hookrightarrow \mathbb P^2$.

A complete description of the monodromy group of an almost generic cover of $\mathbb P^2$
and the local monodromy groups at the points of its branch curve is given by the following
\begin{thm} \label{main} The monodromy group $G_f$ of an  almost generic cover $f:X\to \mathbb P^2$ coincides with
$\mathbb S_{\deg f}$.

The branch curve $B$ of the cover $f$ can have only the singular points of type $A_n$ and the points of $B$ are divided into three types according to the types of singularities of $B$ at these points and properties of the local monodromy groups:
\begin{itemize}
\item[$(i)$]
$p\in B\setminus Sing\, B$ and $G_{f,p}\simeq \mathbb Z_2$ is generated by a transposition;
\item[$(ii)$] $p\in Sing\, B$ is of type $A_{n,2}$, that is {\rm (}by definition{\rm ),} $B$ has the singularity of type $A_{2n-1}$ at $p$, $n\in\mathbb N$, and $G_{f,p}\simeq\mathbb Z_2\times\mathbb Z_2$ is generated by two commuting transpositions;
\item[$(iii)$] $p\in Sing\, B$ is of type $A_{n,3}$, that is {\rm (}by definition{\rm ),} $B$ has the singularity of type $A_{3n-1}$ at $p$, $n\in\mathbb N$, and $G_{f,p}\simeq \mathbb S_3$ is generated by two non-commuting transpositions.
\end{itemize}
\end{thm}

Proof of Theorem \ref{main} is based on the following complete classification of the germs of three-sheeted smooth finite covers.

\begin{thm}\label{T} Let $(U,o')$ and $(V,o)$ be two connected germs of smooth complex-analytic surfaces and $f: (U,o')\to (V,o)$ a finite three-sheeted cover.
Then there are  local coordinates $z,w$ in $(U,o')$ and $u,v$ in $(V,o)$, and a non-negative integer $n\in \mathbb Z_{\geq 0}$ such that $f$ coincides with the cover $f_n:(U,o')\to (V,o)$ given by
$$\begin{array}{l}
u=z, \\
v=w^3-nz^nw, \qquad n\in \mathbb Z_{\geq 0}.
\end{array} $$
\end{thm}

Theorems \ref{main} and \ref{T} are proved in Section 1.

Let $f:X\to \mathbb P^2$ be an almost generic cover branched along a curve $B$.
Below, we use the following notations:
\begin{itemize}
\item $d:=\frac{1}{2}\deg B$;
\item $n_k$, the number of singular points of $B$ of type $A_{2k+1,3}$, $k\in \mathbb Z_{\geq 0}$;
\item $m_k$, the number of singular points of $B$ of type $A_{2k,3}$, $k\in \mathbb N$;
\item $t_k$, the number of singular points of $B$ of type $A_{k,2}$, $k\in \mathbb N$;
\item $\frak{c}= \sum_{k=0}^{\infty}((2k+1)n_k+2km_k),\quad \frak{n}= \sum_{k=1}^{\infty}kt_k,\quad \frak{s}= \sum_{k=1}^{\infty}k(n_k+m_k)$.
\end{itemize}

Note that the numbers $\frak{c}$, $\frak{n}$, and $\frak{s}$ are well-defined, since only for finitely many $k$ the numbers $n_k$, $m_k$, and $t_k$  do not vanish.

In Section 2, we compute the squares of canonical class and the Euler characteristic of the structure sheaf of the covering surface of an almost generic cover $f:X\to\mathbb P^2$ with irreducible branch curve $B\subset \mathbb P^2$  in terms of degree, numbers and singularity types of
$B$. The obtained formulas coincide with similar formulas for generic covers of the plane if we replace $c$ (the number of ordinary cusps) with  $\frak{c}$ and  $n$ (the number of ordinary nodes) with $\frak{n}$ in the formulas for generic covers of the plane in \cite{K1} (compare Claim \ref{R^2}, Propositions \ref{K^2} -- \ref{chi}, Corollary \ref{coro}, and Claim \ref{dual} in Section 2 with Lemmas 4, 6 -- 8, Corollary 2, and the formula for the degree of the dual curve of the branch curve $B$ in \cite{K1}). Therefore, we call $\frak{c}$ the {\it number of pseudo-cusps} and $\frak{n}$ the {\it number of pseudo-nodes} of the branch curve $B$. And in view of Claim \ref{g(B)}, the number $\frak{s}$ is called the {\it superabundance}. Also in Section 2, we investigate the singular points of the Galoisations of almost generic covers of the plane and calculate main invariants of the desingularisations of their covering surfaces.

\section{Proof of Theorems \ref{main} and \ref{T}}
\subsection{Covers $f_n$, $n\geq1$.}\label{subsec1} Consider a finite cover $f_n:\mathbb C^2\to\mathbb C^2$, $n\geq 1$. If we perform the coordinate change
 $(z_1,w_1)=(\sqrt[n]{\frac{n}{3}}z, w)$, $(u_1,v_1)=(\sqrt[n]{\frac{n}{3}}u, v)$, then in new coordinates the cover $f_n$ is given by
\begin{equation} \label{eq00}
\begin{array}{l}
u_1=z_1, \\
v_1=w_1^3-3z_1^nw_1, \qquad n\in \mathbb N.
\end{array} \end{equation}

\begin{claim} \label{cl0} The branch curve $B_n$
of the cover $f_n$, $n\geq 1$, has the singularity of type $A_{3n-1}$ at the point $o=(0,0)$. \end{claim}
\proof The ramification curve $R_n$ of $f_n$  is given by equation $$J(f_n):= \det \left(\begin{array}{cc} 1 & 0 \\ \frac{\partial v_1}{\partial z_1} & \frac{\partial v_1}{\partial w_1}\end{array}\right)=0, $$
i.e., $R_n$ is given by equation
\begin{equation} \label{eqR} w_1^2-z_1^n=0.\end{equation}

Let $n=2k+\delta$, where $\delta= 0$ or $1$ depending on the parity of $n$.

If $n=2k$ is an even number,  then $R_{2k}$ consists of two irreducible components, $R_{2k}=R_+\cup R_-$, given by equations $w=\pm z^{k}$. Therefore the ramification curve $B_{2k}=B_+\cup B_-$, where, due to (\ref{eq00}), $B_+=f_n(R_+)$ and $B_-= f_n(R_-)$ parametrically given by equations
$$ \begin{array}{l}
u_1=z_1, \\
v_1=\mp 2z_1^{3k}
\end{array}$$
and hence, $B_{2k}$ is given by equation
\begin{equation} \label{B2k}
v_1^2-4u_1^{6k}=0
\end{equation}
i.e., in the case when $n=2k$, the point $o$ is the singular point of $B_{n}$ of type $A_{3n-1}$.

If $n=2k+1$ is an odd number, then $R_{n}$ is an irreducible curve  given by equation $w^2-z^{n}=0$. Therefore, $R_{n}$ can be given parametrically by equations $z=t^2$ and $w=t^{2k+1}$ and, due to (\ref{eq00}), the ramification curve $B_{n}$ parametrically given by equations
$$ \begin{array}{l}
u_1=t^2, \\
v_1=-2t^{6k+3}
\end{array}$$
and hence, $B_{n}$ is given by equation
\begin{equation} \label{B2k+1}
v_1^2-4u_1^{6k+3}=0,
\end{equation}
i.e.,  in the case when $n=2k+1$ the point $o$ is the singular point of $B_{n}$ of type $A_{3n-1}$.
 \qed \\

 \subsection{Proof of Theorem \ref{T}.}
Denote $M_1=f^{-1}(L_1)$ and $M_2=f^{-1}(L_2)$ for $L_1=\{u_0=0\}$ and $L_2=\{v_0=0\}$, where $u_0,v_0$ are some local complex-analytic coordinates in $(V,o)$. Then the local intersection number of the curves $M_1$ and $M_2$ at the point $o'$ is equal to $(M_1,M_2)_{o'}=\deg_{o'} f=3$. Therefore either $M_1$ or $M_2$ is a germ of a non-singular curve. Let $M_1$ be non-singular. Then we can choose local coordinates $z_0,w_0$ in $(U,o')$ such that $f^*(u_0)=z_0$ and $f^{*}(v_0)=v_0(z_0,w_0)=\sum_{i=0}^{\infty}a_i(z_0)w_0^i$, where
$a_i(z_0)=\sum_{j=0}^{\infty}a_{i,j}z_0^j\in \mathbb C[[z_0]]$. Performing the coordinates change $v_0 \leftrightarrow v_0-a_0(u_0)$, we can assume that $a_0(z_0)\equiv 0$. In addition, we have $a_{1,0}=a_{2,0}=0$ and can assume that $a_{3,0}=1$, since $(M_1,M_2)_{o'}=3$.

Denote  $R\subset (U,o')$ the ramification curve of the cover $f$ and $B=f(R)\subset (V,o)$ the branch curve. Note that
the curves $R\subset (U,o')$ and $B\subset (V,o)$ depend only on the cover $f$ and do not depend on the choice of coordinates in $(U,o')$ and $(V,o)$. Denote $\mathfrak{m}\subset \mathbb C[[z_0,w_0]]$ the maximal ideal in the ring of power series $\mathbb C[[z_0,w_0]]$.
The curve $R$ is given by equation $$J(f):= \det \left(\begin{array}{cc} 1 & 0 \\ \frac{\partial v_0}{\partial z_0} & \frac{\partial v_0}{\partial w_0}\end{array}\right)=0, $$
i.e., $R$ is given by equation
\begin{equation}\label{eq1} \sum_{i=1}^{\infty}ia_i(z_0)w_0^{i-1}=0.\end{equation}
Let us write equation (\ref{eq1}) in the following form
\begin{equation} \label{eq2} a_{1,1}z_0+a_{1,2}z_0^2+2a_{2,1}z_0w_0 +3w_0^2+ H(z_0,w_0)=0,
\end{equation}
where $H(z_0,w_0)\in \mathfrak{m}^3$. It follows from (\ref{eq2}) that there are three possibilities:
either $R=2R_1$, where $R_1$ is a germ of a smooth curve, or  $R=R_1\cup R_2$, where $R_1$ and $R_2$ are germs of smooth curves, or $R$ is an irreducible germ and $o'$ is its singular point (if $a_{1,1}=0$) of multiplicity two. In the first case the cover $f$ is ramified along $R_1$ with multiplicity three.

\begin{claim}\label{cl-1} If the finite cover $f$ is ramified along $R_1$ with multiplicity three then the branch curve $B$ is smooth and $f$ coincides with the cover $f_0$.
\end{claim}
\proof
It follows from (\ref{eq2}) that the germ $R_1$ is given by equation of the form
$$\sqrt{3}(w_0+\alpha z_0+ H_1(z_0,w_0))=0,$$
where $\alpha\in \mathbb C$ and $H_1(z_0,w_0)\in \mathfrak{m}^2$. In the new system of coordinates $z_1=z_0$, $w_1=w_0+\alpha z_0+ H_1(z_0,w_0)$,
the cover $f$ is given by functions $u=z_1$ and $v_0=v_1(z_1,w_1)$, where
$v_1(z_1,w_1)$ has the following property: $$\frac{\partial v_1}{w_1}=w_1^2(3+H_2(z_1,w_1)),$$
with some $H_2(z_1,w_1)\in \mathfrak{m}$. Therefore $v_0= H_0(z_1)+w_1^3(1+\frac{1}{w_1^3}\int w_1^2H_2(z_1,w_1)dw_1)$, where
$\frac{1}{w_1^3}\int w_1^2H_2(z_1,w_1)dw_1\in \mathfrak{m}$, and it is easy to see that in the new systems of coordinates $$z=z_1,\,\, w=w_1\sqrt[3]{1+\frac{1}{w_1^3}\int w_1^2H_2(z_1,w_1)dw_1},\quad  \text{and}\quad u=u_0,\,\,v=v_0-H_0(u_0)$$ the cover $f$ coincides with $f_0$.  Note that the branch curve $B$ of the cover $f_0$ is smooth and it is given by equation $v=0$. \qed\\

Now, we assume that the cover $f$ is ramified along $R$ with multiplicity two.
\begin{claim} \label{cla} The restriction $f_{\mid R}:R\to B$ of the cover $f$ to the ramification locus $R$ is one-to-one mapping.
\end{claim}
\proof Obviously, $3=\deg f \geq 2\deg f_{\mid R}$. Therefore $\deg f_{\mid R}=1$. \qed

\begin{claim} \label{cl1}  Let the branch curve $B\subset (V,o)$ of the cover $f$ have a singularity of type $A_m$ at the point $o$. Then  $m=3n-1$ for some $n\in\mathbb N$ and $f$ coincides with the cover $f_n$.
\end{claim}
\proof The cover $f:(U,o')\to (V,o)$ defines a homomorphism $$f_*:\pi_1^{loc}(B,o)=\pi_1(V\setminus B,p)\to \mathbb S_3,$$ where $\mathbb S_3$ is the symmetric group acting on the fibre $f^{-1}(p)$. Note that the epimorphism $f_*$ is defined uniquely only if we fix a numbering of the points of $f^{-1}(p)$ and in general case it defines uniquely up to inner automorphism of $\mathbb S_3$.

It is well known (see, for example, \cite{K-U}) that if $m=2k-\delta$, where $\delta= 0$ or $1$, then the group $\pi_1^{loc}(B,o)$ is generated by two so called geometric generators $\gamma_1$ and $\gamma_2$ such that $\pi_1^{loc}(B,o)$ has the following presentation:
\begin{equation}\label{pres}\pi_1^{loc}(B,o)=\langle \gamma_1, \gamma_2\mid  (\gamma_1\gamma_2)^k\gamma_1^{1-\delta}=(\gamma_2\gamma_1)^k\gamma_2^{1-\delta}\rangle .\end{equation}
Denote $\tau_i:=f_*(\gamma_i)\in \mathbb S_3$, $i=1,2$. The branch curve $B$ is singular. Therefore, by Claim \ref{cl-1}, $\tau_i$ are not cycles of length three, i.e.,  $\tau_1$ and $\tau_2$ are transpositions generating the group $\mathbb S_3$, since $(U,o')$ is a germ of an irreducible surface. Without loss of generality, we can assume that $\tau_1=(1,3)$ and $\tau_2=(2,3)$, i.e., up to conjugation  there is the unique epimorphism from $\pi_1^{loc}(B,o)$ to $\mathbb S_3$.
We have $\tau_1\tau_2=(1,2,3)$, $\tau_2\tau_1=(1,3,2)$ and it follows from (\ref{pres}) that
$$(1,2,3)^k(1,3)^{1-\delta}=(1,3,2)^k(2,3)^{1-\delta}.$$
If $\delta=1$ then $k=3k_1$, i.e., $m=3(2k_1)-1$ and if $\delta=0$ then $k=3k_1-2$, i.e.,
$m=2(3k_1-2)=3(2k_1-1)-1$, where $k_1\geq 1$. By  Grauert - Remmert - Riemann - Stein Theorem (\cite{St}, \cite{G-R}), the cover $f:(U,o')\to (V,o)$ is uniquely defined  by epimorphism $f_*:\pi_1^{loc}(B,o)=\pi_1(V\setminus B,p)\to \mathbb S_3$. Now, to complete the proof of Claim \ref{cl1}, it suffices to apply Claim \ref{cl0}. \qed \\

It follows from Claim \ref{cl1} that to prove Theorem \ref{T}, it suffices to show that in the second and third cases the branch curve $B$ has the singularity of type $A_m$ for some $m\geq 1$. Therefore it suffices to show (see, for example, \cite{Ar}) that the multiplicity of the singular point $o$ of $B$ is equal to $2$.

It follows from (\ref{eq2}) that in the second case $a_{1,1}=0$ and the germ $R=R_1\cup R_2$ is the union of two curves smooth at $o'$. In addition, it is easy to see that  equations of $R_i$, $i=1,2$, have the following form
$$ w_0+\alpha_iz_0+ H_i(z_0,w_0)=0, $$
where $H_i(z_0,w_0)\in \mathfrak{m}^2$ and $\alpha_i\in \mathbb C$. Therefore the function $z_0$ is a local parameter at $o'$ for the germs $R_i$, $i=1,2$, i.e., $B_i=f(R_i)\subset V$ are germs of  smooth curves at the point $o$, since $u_0=z_0$, and, hence the multiplicity of the singular point $o$ of $B$ is equal to $2$. As a result, we obtain that $B$ has a singularity of type $A_{2m-1}$ for some $m$ equals to the intersection number $(B_1,B_2)_{o}$.

In the third case, denote by $\nu:\widetilde R\to R$ the resolution of singular point $o$ of $R$ (if $a_{1,1}\neq 0$ then $\nu=id$), $p=\nu^{-1}(o')$, and denote by $\frak{m}_{\widetilde R,p}$ the maximal ideal of the ring of holomorphic functions at $p$ on
$\widetilde R$.

It follows from (\ref{eq2}) that  $(L_1,R)_{o'}=2$. Therefore $\nu^*(u_0)\in \mathfrak{m}_{\widetilde R,p}^2\setminus \mathfrak{m}_{\widetilde R,p}^3$ and the function $t=\sqrt{\nu^*(u_0)}$ is a local parameter on $\widetilde R$ at $p$. Let $\nu^*(w_0)=w(t)\in \frak{m}_{\widetilde R,p}$. Then $R$ is given parametrically  by
$$ z_0=t^2, \qquad w_0=w(t)$$ and therefore $B$ is given by
\begin{equation} \label{param}
u_0=t^2, \qquad v_0=\sum_{i=0}^{\infty} \sum_{j=0}^{\infty} a_{i,j}t^{2j}w(t)^i, \end{equation}
where $v_0(t)\in \frak{m}_{\widetilde R,p}^3$, since $a_{1,0}=a_{2,0}=0$. It follows from (\ref{param}) that the multiplicity of $B$
the point $o$ is equal to $2$ and therefore the singularity type of $B$ at $o$ is $A_{2m}$ for some $m\in \mathbb N$. \qed

\subsection{Proof of Theorem \ref{main}.}
Proof of the following Claim will be omitted, since it is similar and much more simpler than the proof of Theorem \ref{T}.
\begin{claim}\label{2cl} Let $(U,o')$ and $(V,o)$ be two germs of smooth complex-analytic surfaces and $f: (U,o')\to (V,o)$ a finite two-sheeted cover.
Then there are  local coordinates $z,w$ in $(U,o)$ and $u,v$ in $(V,o)$ such that $f$ is given by equations
$u=z$, $v=w^2$.

The ramification locus $R\subset (U,o')$ and the branch curve $B\subset (V,o)$ of the cover $f$ are smooth curves and they are given, respectively, by equations $w=0$ and $v=0$.
\end{claim}

Let $V\subset \mathbb P^2$ be a sufficiently small neighbourhood of a point $q$ of the branch curve $B$ of an almost generic cover $f$. Let $\deg f=N$. It follows from property $(G_1)$ that there are three possibilities:
\begin{itemize}
\item[$(1)$] $f^{-1}(V)$ is a disjoint union of $N-1$ neighbourhoods $U_1,\dots, U_{N-1}$ such that $f:U_1\to V$ is a two-sheeted cover and $f:U_i\to V$ are biholomorphic maps for $i=2,\dots, N-1$;
\item[$(2)$] $f^{-1}(V)$ is a disjoint union of $N-2$ neighbourhoods $U_1,\dots, U_{N-2}$ such that $f:U_1\to V$ and $f:U_2\to V$ are  two-sheeted covers and $f:U_i\to V$ are biholomorphic maps for $i=3,\dots, N-2$;
\item[$(3)$] $f^{-1}(V)$ is a disjoint union of $N-2$ neighbourhoods $U_1,\dots, U_{N-1}$ such that $f:U_1\to V$ is a three-sheeted cover and $f:U_i\to V$ are biholomorphic maps for $i=2,\dots, N-1$.
\end{itemize}

\begin{claim} \label{birat} Let $B_1$ be an irreducible component of the branch curve $B\subset \mathbb P^2$ of an almost generic cover $f$ and $R_1=f^{-1}(B_1)\cap R$. Then the restriction $f_{\mid R_1}:R_1\to B_1=f(R_1)$ of $f$ to $R_1$
is a birational morphism.
\end{claim}
\proof It follows from property $(G_1)$ that $\deg f_{\mid R_1}\leq 2$. Assume that $\deg f_{\mid R_1}= 2$. Then it follows from Claims \ref{2cl} and \ref{cla} that $B_1$ is a smooth curve and $f_{\mid R_1}:R_1\to B_1$ is unramified two-sheeted cover. Applying  property $(G_1)$, we obtain that $B=B_1$, since any two curves in $\mathbb P^2$ have non-empty intersection. Therefore, by Zariski Theorem, the group $\pi_1(\mathbb P^2\setminus B)$ is cyclic and hence the monodromy group $G_f\simeq \mathbb Z_2$ generated in $\mathbb S_{\deg f}$ by product of two commuting transpositions. Therefore $G_f$ can not act transitively on the set $f^{-1}(p)$ 
which contradicts the irreducibility of $X$.  \qed \\

The surface $X$ is irreducible. Therefore $X\setminus f^{-1}(B)$ is connected and hence, $G_f$ acts transitively on the fibre $f^{-1}(p)$ over the base point $p$ of the group $\pi_1(X\setminus f^{-1}(B),p)$. It follows from property $(G_2)$ and Claim \ref{birat} that the monodromy group $G_f\subset \mathbb S_{\deg f}$ of an almost generic cover $f:X\to \mathbb P^2$ is generated by transpositions. Therefore the monodromy group $G_f$ coincides with $\mathbb S_{\deg f}$.

It follows from Claim  \ref{2cl} that in case $(1)$ locally at the point $q$, the branch curve $B$ is smooth (the germ $(B,q)$ is the branch curve of the restriction of $f$ to the neighbourhood $U_1$).
By Claim \ref{2cl}, the local fundamental group $\pi_1^{loc}(B,q)$ is generated by circuit $\gamma$ around $B$ and the local monodromy group $G_{f,q}$ is generated by transposition $\tau=f_*(\gamma)$.

It follows from Claim \ref{2cl} that in case $(2)$ locally at the point $q$, the branch curve $B$ consists of  two smooth branches (denote them by $B_1$ and $B_2$), $B_1$ is the branch curve of the restriction of $f$ to the neighbourhood $U_1$ and $B_2$ is the branch curve of the restriction of $f$ to the neighbourhood $U_2$. Note that by Claim  \ref{birat}, we have $B_1\neq B_2$. Therefore $B$ at the point $q$ has the singularity of type $A_{2k-1}$, where $k=(B_1,B_2)_q$. The local fundamental group $\pi_1^{loc}(B,q)$ is generated by circuits $\gamma_1$ around $B_1$ and $\gamma_2$ around $B_2$ and the local monodromy group $G_{f,q}$ is generated by commuting transpositions $\tau_1=f_*(\gamma_1)$ and $\tau_2=f_*(\gamma_2)$. Therefore $q$ is a point of type $A_{k,2}$.

Using Theorem \ref{T} and Claim \ref{cl0}, it is easy to see that in case $(3)$ the point $q$ is of type $A_{k,3}$ for some $k\in \mathbb N$.

\section{Relations between invariants of covering surfaces of almost generic covers and invariants of their branch curves}
\subsection{Invariants of the branch and ramification curves.}
Let $f:X\to\mathbb P^2$ be an almost generic cover. We will assume that the branch curve $B\subset\mathbb P^2$ of $f$ is irreducible and let $R\subset X$ be the ramification locus of $f$.
\begin{claim} \label{degf} The degree of $B$ is  an even number, $\deg B=2d$, $d\in\mathbb N$.
\end{claim}
\proof Let $L$ be a line in $\mathbb P^2$ generic with respect to $B$ and $M=f^{-1}(L)$. By property $(G_2)$, the restriction  $f_{\mid M}:M\to L$ of $f$ to $M$ is a generic cover branched over the common points of $L$ and $B$ and by Hurwitz formula,
a generic cover of $\mathbb P^1$ is branched over even number of points. Therefore $\deg B=(L,B)_{\mathbb P^2}$ is an even number. \qed \\

The following claim is well known.
\begin{claim}\label{delta} If $s\in Sing\, B$ is of type $A_{k,2}$ then its $\delta$-invariant is equal to $k$; if $s$ is of type $A_{2k+1,3}$, $k\geq 0$, then its $\delta$-invariant is equal to $3k+1$;  if $s$ is of type $A_{2k,3}$, $k\geq 1$, then its $\delta$-invariant is equal to $3k$.
\end{claim}

\begin{claim} \label{g(B)} The geometric genus $g(B)$ of $B$ is equal to
\begin{equation}\label{gen} g(B) =(2d-1)(d-1)-\frak{c}-\frak{n}- \frak{s}.
\end{equation}
\end{claim}
\proof By Theorem \ref{main}, the singular points of the curve $B$ are of the types $A_{n,2}$ and $A_{n,3}$, $n\in \mathbb N$. Therefore Claim \ref{g(B)} follows from Claim \ref{delta}. \qed
\begin{claim} \label{dual} The degree $\hat d=\deg \hat B$ of the dual curve $\hat B$ of $B$ is equal to
$$ \hat d= 2d(2d-1)- 3\frak{c}-2\frak{n}. 
$$
\end{claim}
\proof If $s\in Sing\, B$ is of type $A_{k,2}$ then its number of virtual cusps vanishes and the number of its virtual nodes  is equal to $k$; if $s$ is of type $A_{2k+1,3}$, $k\geq 0$, then its number of virtual cusps is equal to $1$ and its number of  virtual nodes  is equal to  $3k$;  if $s$ is of type $A_{2k,3}$, $k\geq 1$, then its number of virtual cusps vanishes and its number of  virtual nodes  is equal to  $3k$. Now, Claim \ref{dual} follows from generalized Pl$\ddot{\rm{u}}$kker's formula (see \cite{Ku}). \qed

\begin{claim} \label{R^2}
The self-intersection number $(R^2)_X$ of $R$ is positive and it is equal to
\begin{equation}\label{(R,R)} (R^2)_X= 2d^2-\frak{c}-\frak{n}.
\end{equation}
\end{claim}
\proof
It follows from Theorems \ref{main} and \ref{T}, and equation (\ref{eqR}) that $s$ is a singular point of the curve $R$  (and its singular type is $A_{n-1}$) iff $f(s)$ is a point of  type $A_{n,3}$, $n\geq 1$. Therefore
\begin{equation} \label{g(R)} 2(g(R)-1)= (K_X+R,R)_X - 2\sum_{k=1}^{\infty} k(n_k+m_k). \end{equation}
By Claim \ref{birat}, we have $g(B)=g(R)$. In addition, $K_X=f^*(K_{\mathbb P^2})+R$ and $$(f^*(K_{\mathbb P^2}),R)_X=-3\deg B=-6d.$$ Therefore, equality (\ref{(R,R)}) follows from (\ref{gen}) and (\ref{g(R)}) and it follows from Claim \ref{dual} that $(R^2)_X>0$, since $\hat d>0$. \qed

\subsection{Invariants of the covering surfaces.}
Let $N$ be the degree of an almost generic cover $f:X\to\mathbb P^2$.
\begin{prop} \label{K^2}
The self-intersection number $K^2_X$ is equal to
\begin{equation} \label{kan}
K_X^2=9N +2(d^2-6d)-\frak{c}-\frak{n}.
\end{equation}
\end{prop}
\proof We have $K_X=f^*(K_{\mathbb P^2})+R$. Therefore, by Claim \ref{R^2},
$$\begin{array}{ll}  K_X^2= & (f^*(K_{\mathbb P^2}),f^*(K_{\mathbb P^2}))_X+ 2(f^*(K_{\mathbb P^2}),R)_X +R^2= \\ &
\displaystyle 9N-12d + 2d^2-\frak{c}-\frak{n}.
\qquad \qed \end{array} $$ \vspace{0.1cm}

Denote by $e(M)=\sum_{i}(-1)^i\dim H^i(M,\mathbb Q)$ the topological Euler characteristic of a topological space $M$.

\begin{prop} \label{e(X)}
The topological Euler characteristic $e(X)$ is equal to
$$ e(X)=3N+ 2d(2d-3)- 3\frak{c} -2\frak{n}.
$$
\end{prop}
\proof We have
$e(X)=N(e(\mathbb P^2) - e(B))+(N-1)e(B\setminus Sing\, B)+(N-2)e(Sing\, B)$. Therefore Proposition \ref{e(X)} follows from (\ref{gen}) and equalities
\begin{equation} \label{e(B)} e(B)  =(2-2g(B))-\sum_{k=1}^{\infty}(m_k+t_k), \end{equation}
\begin{equation} \label{e(B)-} e(B\setminus Sing\, B)  =(2-2g(B))-n_0-\sum_{k=1}^{\infty}(n_k+2m_k+2t_k), \end{equation}
\begin{equation} e(Sing\, B)  = n_0+\sum_{k=1}^{\infty}(n_k+m_k+t_k).
\end{equation}

\begin{prop} \label{chi} The Euler characteristic $\chi(\mathcal O_X)$ of the structure sheaf $\mathcal O_X$ equals
$$\chi(\mathcal O_X)=N+\frac{d(d-3)}{2}-\frac{\frak{c}}{3}- \frac{\frak{n}}{4}.$$
\end{prop}
\proof Proposition \ref{chi} follows from Noether's formula $K_X^2+e(X)=12\chi(\mathcal O_X)$. \qed
\begin{cor} \label{coro} We have
$$ \frak{c}\equiv 0\,\, (\text{ \rm
mod}\,\, 3), \quad \frak{n}\equiv 0\,\, (\text{\rm mod}\,\, 4). $$
\end{cor}

\begin{prop}\label{Hodge} The following inequality holds
$$N \leq \frac{4d^2}{2d^2- \frak{c}-\frak{n}
}\, .$$
\end{prop}
\proof
Let $M=f^{-1}(L)$, where $L$ is a line in$\mathbb P^2$. We have $(M^2)_X=N$ and $(M,R)_X=\deg f= 2d$.
Applying  Hodge Index Theorem to $M$ and $R$ and applying Claim \ref{R^2}, we
obtain  inequality
\begin{equation}
N( 2d^2- \frak{c}-\frak{n}
) -4d^2\leq 0.\qquad \qed    \label{h}
\end{equation}

\begin{claim} \label{N>} If $N\geq 6$ then $\mathfrak{n}>0$. \end{claim}
\proof Transforming the left side of inequality (\ref{h}), we obtain the inequality
$$
\frac{1}{2}[(N-6)(2d(2d-1)+2d- 3\frak{c}-2\frak{n}+\frak{c})
+4(2d(2d-1)- 3\frak{c}-2\frak{n}+2d) -4\frak{n}]\leq 0
$$
and applying Claim \ref{dual}, we have
$$(N-6)(\hat d+2d+\frak{c})
+4(\hat d+2d)\leq 4\frak{n}.
$$
Now, Claim \ref{N>} follows from inequalities $\hat d>0$ and $d>0$.\qed

\subsection{Galoisations of almost generic covers.}
The Cayley imbedding $c: G_f=\mathbb S_N\hookrightarrow \mathbb S_{N!}$, defined by the action of $\mathbb S_N$ on itself by multiplication from the right side, defines the Galois finite cover $g:Y\to \mathbb P^2$ branched along the curve $B$,
$\deg g=N!$. For a point $p\in \mathbb P^2$
the fibre $g^{-1}(p)$ is supported on $M_{f,p}=\frac{N!}{|G_{f,p}|}$ distinct points, where $|G_{f,p}|$ is the order of the local monodromy group $G_{f,p}=G_{g,p}$, and if $V\subset \mathbb P^2$ is a sufficiently small complex-analytic neighbourhood of the point $p$, then $g^{-1}(V)= \bigsqcup_{i=1}^{M_{f,p}}W_i$ is a disjoint union of $M_{f,p}$ complex-analytic normal varieties $W_i$ biholomorphic to each other and the restriction $g_{\mid W_i}:W_i\to V$ to each $W_i$ is the Galois cover with the Galois group isomorphic to $G_{f,p}$.

If $p\in B\setminus Sing\,B$ then for each $i=1,\dots, M_{f,p}=\frac{N!}{2}$ the cover $g_{\mid W_i}:W_i\to V$ is a two-sheeted cover branched over the non-singular curve $B\cap V$.

If $p\in Sing\,B$ and $B$ has the singularity of type $A_{n,2}$ at $p$,  then for each $i=1,\dots, M_{f,p}=\frac{N!}{4}$ the cover $g_{\mid W_i}:W_i\to V$ is the Galois cover with Galois group $\mathbb Z_2\times \mathbb Z_2$ branched along $B\cap V$ with multiplicity two.

\begin{claim} \label{Gal1} If $B$ has the singularity of type $A_{n,2}$ at $p$,  then for $n\geq 2$ the point $g_{|W_i}^{-1}(p)$ is a singular point of $W_i$ of type $A_{n-1}$ and it is smooth point if $n=1$. \end{claim}
\proof Let $B_1$ and $B_2$ be two irreducible branches of the curve $B\cap V$. Without loss of generality, we can assume that $u=0$ is the equation of $B_1$ and $u-v^n=0$ is the equation of $B_2$. Then $W_i$ is biholomorphic to a neighbourhood of the point $o=(0,0,0,0)$ in the surface given  in $\mathbb C^4$ by equations $z^2=u$ and $w^2=u-v^n$. Therefore, it is biholomorphic  to a neighbourhood of the point $o'=(0,0,0)$ in the  surface given in $\mathbb C^3$ by equation $w^2=z^2-v^n$. \qed\\

If $p\in Sing\,B$ and $B$ has the singularity of type $A_{n,3}$ at $p$,  then for each $i=1,\dots, M_{f,p}=\frac{N!}{6}$ the cover $g_{\mid W_i}:W_i\to V$ is the Galois cover with Galois group $\mathbb S_3$ branched along $B\cap V$ with multiplicity two.

\begin{claim} \label{Gal2} If $B$ has the singularity of type $A_{n,3}$ at $p$,  then for $n\geq 2$ the point $g_{\mid W_i}^{-1}(p)$ is a singular point of $W_i$ of type $A_{n-1}$ and it is smooth point if $n=1$. \end{claim}
\proof The local fundamental group $\pi_1^{loc}(B,p)$ is generated by two circuits $\gamma_1$ and $\gamma_2$ around the curve $B$ (see subsection \ref{subsec1}) and the Galois cover $g_{\mid W_i}:W_i\to V$ is defined by epimorphism $g_*:\pi_1^{loc}(B,p)\to \mathbb S_3$ sending $\gamma_i$, $i=1,2$, to transpositions.
The cover $g_{\mid W_i}:W_i\to V$ can be decomposed into a composition $g_{\mid W_i}=f_n\circ h_n$, where $h_n:W_i\to W_i/ \langle \tau \rangle\simeq U\subset\mathbb C^2$ is the factor-map under the action of a subgroup $\langle \tau \rangle$ of the group $\mathbb S_3$ generated by a transposition $\tau\in\mathbb S_3$ and by Theorem \ref{T}, $f_n:U\to V$ is the restriction of the cover $f_n:\mathbb C^2\to \mathbb C^2$ (see subsection \ref{subsec1}) to some neighbourhood $U\subset \mathbb C^2$.

It is easy to see that $h_n:W_i\to U\subset\mathbb C^2$ is a two-sheeted cover branched along the curve $C_n\cap U$, where $C_n$ is the complement to the ramification divisor $R_n$ of the cover $f_n:\mathbb C^2\to\mathbb C^2$ in the total inverse image of $B_n$,  $f_n^*(B_n)=2R_n+C_n$. Let $F_n(z_1,w_1)=0$ be an equation of $C_n$. Then it follows from (\ref{eqR}) -- (\ref{B2k+1}) that
$$(w_1^3-3z_1^nw_1)^2-4z_1^{3n}=(w_1^2-z_1^n)^2F_n(z_1,w_1)=0. $$
Therefore we have
$ F_n(z_1,w_1)=w_1^2-4z_1^n$ and hence the surface $W_i$ is isomorphic to a hypersurface in $\mathbb C^1\times U$ given by
$ y^2=w_1^2-4z_1^n$. \qed \\

It follows from Claims \ref{Gal1} and \ref{Gal2} that the covering surface $Y$ of the Galoisation of an almost generic cover $f:X\to\mathbb P^2$ has
$$ S=N!(\frac{1}{6}\sum_{k=1}^{\infty}(n_k+m_k)+\frac{1}{4}\sum_{k=2}^{\infty}t_k)$$
singular points.

Let $\nu: Z\to Y$ be the minimal resolution of singular points of the covering surface $Y$.
It is well known that the inverse image $\nu^{-1}(s)$ of a point $s\in Sing\, Y$ of singularity type $A_k$ is the chain $E_1\cup\dots\cup E_k$ of $(-2)$-curves $E_i$. Therefore, by Claims \ref{Gal1} and \ref{Gal2}, the number of $(-2)$-curves contracted by $\nu$ is
$$M= \frac{N!}{6}\sum_{k=1}^{\infty}(2kn_k+(2k-1)m_k)+ \frac{N!}{4}\sum_{k=2}^{\infty}(2k-1)t_k$$
and if we denote $\widetilde g:=g\circ \nu:Z\to \mathbb P^2$, then
\begin{equation} \label{e(S)}
e(\widetilde g^{-1}(Sing\, B))=\frac{N!}{6}\sum_{k=0}^{\infty}((2k+1)n_k+2km_k)+ \frac{N!}{4}(t_1+\sum_{k=2}^{\infty}2kt_k).\end{equation}

\begin{prop} \label{K^2g} The canonical class $K_Z$ of the surface $Z$ is equal to $$K_Z=(d-3)\widetilde g^*(L),$$ where $L$ is a line in $\mathbb P^2$, and its self-intersection number $K^2_Z$ is equal to
\begin{equation} \label{kang}
K_Z^2=(d-3)^2N!.
\end{equation}
\end{prop}
\proof Let $\nu^{-1}(Sing\, Y)=\cup_{i=1}^ME_i$. We have
$$K_Z=\widetilde g^*(K_{\mathbb P^2})+\frac{1}{2}\widetilde g^*(B)+\sum_{j=1}^M\alpha_jE_j,$$
since $\widetilde g$ is ramified over $B$ with multiplicity two. In addition, we have $(K_Z,E_i)_Z=0$ for each $i$, since $E_i$ are rational curves with self-intersection number $(E_i^2)_Z=-2$, and $(\widetilde g^*(K_{\mathbb P^2}),E_i)_Z=(\widetilde g^*(B),E_i)_Z=0$, since $\widetilde g(E_i)$ are points. Therefore $(\sum_{j=1}^M\alpha_jE_j,E_i)_Z=0$ for each $i$ and hence $\sum_{i=1}^M\alpha_iE_i=0$, since the intersection matrix $E=((E_j,E_i)_Z)$ is negatively defined. Now, to complete the proof of Proposition \ref{K^2g}, notice that $K_{\mathbb P^2}=-3L$, the divisor $B$ is equivalent to $2dL$, and $\deg \widetilde g= N!$. \qed

\begin{prop} \label{e(Z)}
The topological Euler characteristic $e(Z)$ is equal to
$$ e(Z)=N![3+ d(2d-3)- \displaystyle \frac{1}{6}\sum_{k=0}^{\infty}((19k+5)n_k+ 
16km_k)-\frac{3}{4}t_1-\frac{1}{2}\sum_{k=2}^{\infty}kt_k].$$
\end{prop}
\proof We have
$$e(Z)=N!(e(\mathbb P^2) - e(B))+\frac{N!}{2}e(B\setminus Sing\, B)+e(\widetilde g^{-1}(Sing\, B))$$
and Proposition \ref{e(Z)} follows from (\ref{gen}), (\ref{e(B)}), (\ref{e(B)-}), and (\ref{e(S)}). \qed \\

\ifx\undefined\bysame
\newcommand{\bysame}{\leavevmode\hbox to3em{\hrulefill}\,}
\fi


\begin{thebibliography}{99}
\bibitem{Ar} V.I. Arnol'd: {\it Normal forms for functions near degenerate critical points, the Weyl groups of $A_k$, $D_k$, $E_k$
 and Lagrangian singularities,} Funct. Anal. Appl., {\bf 6:4} (1972), 254 -- 272.


\bibitem{G-R} {H. Grauert, R. Remmert:}
{\it Komplexe R$\ddot{\text{a}}$ume,} Math. Ann., {\bf 136}(1958), 245 -- 318.

\bibitem
{C-F} {C. Ciliberto, F. Flamini:} {\it On the branch curve of a generic projection of a surface to a plane,} arXiv:0811.0467

\bibitem{Ku} Vik. S. Kulikov: {\it A Remark on Classical Pl\"ucker's formulae}, Ann. Fac. Sci. Toulouse. Math., 25:5 (2016),
959 -- 967.

\bibitem {K1}  {Vik.S. Kulikov:} {\it On Chisini's conjecture,} Izv. Math., {\bf 63:6} (1999), 1139 -- 1170.

\bibitem{K2} Vik.S. Kulikov: {\it On Chisini's conjecture. II,} Izv. Math., {\bf 72:5} (2008), 901 -- 913.

\bibitem{K3} {Vik.S. Kulikov:} {\it Generalized Chisini's Conjecture,} Proc. Steklov Inst. Math., {\bf 241} (2003), 110 -- 119.

\bibitem{K-U} {Vik.S. Kulikov:} {\it Hurwitz curves,} Russian Math. Surveys, {\bf 62:6} (2007), 1043 -- 1119.

\bibitem{N} S.Yu. Nemirovski: {\it Kulikov’s theorem on the Chisini conjecture,} Izv. Math., {\bf 65:1}
(2001), 71 -- 74.

\bibitem{No} M. Nori: {\it Zariski’s conjecture and related problems}, Ann. Sci. $\rm \acute{E}$cole Norm. Sup.
Ser. {\bf 4, 16} (1983), 305 -- 344.

\bibitem{St}{K. Stein:} {\it Analytische Zerlegungen komplexer R\"aume,}
Math. Ann. \textbf{132} (1956), 63--93.

\end{thebibliography}
\end{document}